\documentclass{article}
\usepackage[utf8]{inputenc}

\usepackage{hyperref}
\hypersetup{
    colorlinks=true,
    linkcolor=blue,
    filecolor=magenta,      
    urlcolor=blue,
    citecolor=black
}

\title{Convex Relaxations \\ in Power System Optimization: \\ A Brief Introduction}

\author{Carleton Coffrin$^1$ \and Line Roald$^2$}
\date{%
    $^1$Los Alamos National Laboratory\\%
    $^2$University of Wisconsin-Madison\\[2ex]%
    \today
}

\begin{document}

\maketitle

\section*{Overview \\
\small{(\href{https://www.youtube.com/watch?v=gB43TmcoUpA&list=PLeuOzWTGxj2ZZ_XUutDwNFvNfSWwWCgR5}{Link to the video series})}}

Convex relaxations of the AC power flow equations have attracted significant interest in the power systems research community in recent years.  The following collection of video lectures provides a brief introduction to the mathematics of AC power systems, continuous nonlinear optimization, and relaxations of the power flow equations. The aim of the videos is to provide the high level ideas of convex relaxations and their applications in power system optimization, and could be used as a starting point for researchers who want to study, use or develop new convex relaxations for use in their own research. The videos do not aim to provide an in-depth tutorial about specific convex relaxations, but rather focus on ideas that are common to all convex relaxations of the AC optimal power flow problem.

The videos assume minimal prerequisite knowledge and are designed to be accessible to a variety of disciplines (e.g. power system engineers, industrial engineers and  computer scientists). 
They are designed to be modular, such that interested listeners can either watch the series as a whole, or choose only the subset that matches their interest. The following sections provide a brief overview of the topics and references covered in each video.

\subsection*{Preliminaries \\ 
\small{(\href{https://youtu.be/gB43TmcoUpA}{Link to the video})}}
This video provides the introduction to the video series. It begins with the motivation, goals, and scope of the series and then introduces some mathematical foundations of power systems and optimization. \cite{glover2011power,9780070359581,grainger1994power}

\subsection*{AC Power Flow \\
\small{(\href{https://youtu.be/Yr1kbZxJ3Eo}{Link to the video})}} 
This video introduces the basics of power networks, and presents a stylized version of the AC Power Flow problem.  This simple power flow model is used to provide an intuition for power losses on lines and describe how power flows in meshed networks with cycles.

\subsection*{AC Optimal Power Flow \\
\small{(\href{https://youtu.be/_ILr-0MkRug}{Link to the video})}}
Further extending the AC Power Flow problem, this video introduces the AC Optimal Power Flow problem \cite{ac_opf_origin}. While these lectures focus on the economic dispatch problem, where generator schedules are optimized to minimize generation cost while satisfying network constraints, the broad applicability of the core AC power constraints to a wide variety of power system decision problems is also discussed.

\subsection*{Computational Hardness and the Value of Convexity \\ \small{(\href{https://youtu.be/uvhZXuPCmk0}{Link to the video})}}
It is commonly acknowledged that the AC Optimal Power Flow problem is ``hard''. This video provides a brief introduction to the concept of computational complexity \cite{sipser2013introduction}, a theoretical field which provides a scientifically rigorous  definition of what ``hard'' means. The \emph{NP-hardness} of AC Optimal Power Flow \cite{7063278,1512.07315} is discussed and used to motivate the value of convexity in power system optimization \cite{6815671,6345272}.

%In preparation for a discussion of optimization methods and their implications for the AC Optimal Power Flow problem, this video provides a brief introduction to the concept of computational complexity \cite{sipser2013introduction}, discusses the NP-hardness of AC Optimal Power Flow \cite{7063278,1512.07315} and motivates the value of convexity in power system optimization \cite{6815671,6345272}.

\subsection*{Solution Methods for AC Optimal Power Flow \\ \small{(\href{https://youtu.be/tZ4Gw-572Mw}{Link to the video})}}
This video discusses a variety of traditional solution methods for the AC Optimal Power Flow problem, including local optimization methods, global optimization methods, and linear approximations \cite{LPAC_ijoc,6407493}.  Additionally, the existence of local optimal solutions in the AC Optimal Power Flow problem is briefly discussed \cite{6581918}.

\subsection*{Convex Relaxations of Non-Linear Optimization Problems\\ \small{(\href{https://youtu.be/5A1ih3Z5YbA}{Link to the video})}}
This video introduces the concept of the \emph{feasible set} of an optimization problem \cite{932273,7879340}, as well as convex relaxations and convex restrictions of that feasible set.  Combining the true feasible set and a relaxation, the concept of computing \emph{optimality gaps} and how relaxations can be used to \emph{prove infeasiblity} are discussed.  The presentation concludes with a discussion of what a \emph{good} convex relaxation is and a review of different methods for strengthening relaxations.

\subsection*{Convex Relaxations of AC Optimal Power Flow \\ \small{(\href{https://youtu.be/ONJYq8uPOU4}{Link to the video})}}
Combining results from the previous videos, this video revisits the AC Optimal Power Flow problem and illustrates how a simple convex relaxation \cite{1664986} can be developed for this non-convex problem.  This is followed by a brief overview of different relaxations \cite{forthcoming}
and a discussion of the challenges that convex relaxations face in accurately representing power losses on lines and power flow around cycles in a network \cite{6345272}.  Special cases when a convex relaxation also provides a feasible solution to AC Optimal Power Flow are discussed \cite{5971792,6815671,6822653,7447798}, as well as more challenging cases where simple convex relaxations may fail \cite{6120344,qc_opf_tps,nesta,pglib_opf}.

\subsection*{Tips for Using Convex Relaxations of AC Optimal Power Flow \\ \small{(\href{https://youtu.be/yTbOK6i6_SY}{Link to the video})}}
This final video concludes the series with a collection of best practices for using and developing convex relaxations. The presentation covers several methods for using convex relaxations as building blocks in more complex algorithms, and presents some common caveats and situations where convex relaxations should be used with caution. It also provides advice on how to test whether your convex relaxation is better than existing relaxations.

\subsection*{Extra Materials}

These extra lectures provide examples of how convex relaxations can be leveraged as a building block to develop novel power system analysis tools.

\subsubsection*{Optimization-Based Bound Tightening via Convex Relaxation \\ \small{(\href{https://youtu.be/63rE-kI4xAs}{Link to the video})}}
This lecture provides a brief overview of how convex relaxations can be used for bound tightening in optimal power flow problems, as presented in \cite{cp_qc_fp}.  The presentation illustrates that bound tightening can greatly improve the optimality gaps of convex relaxations and that these tight bounds can provide valuable insights into operational flexibility of a given power network, especially for a fixed operating point.

\subsubsection*{Towards AC Optimal Power Flow with Robust Feasibility Guarantees \\ \small{(\href{https://youtu.be/wpd2JItVibk}{Link to the video})}}
This lecture demonstrates how convex relaxations can be utilized to guarantee robust constraint feasibility for a stochastic variant of the AC Optimal Power Flow problem, as presented in \cite{molzahn2018}. The method uses convex relaxations to provide conservative bounds on the impact of uncertainty. To avoid overly conservative results, a \emph{combination} of \emph{multiple} convex relaxations and bound tightening is utilized.

\section*{Acknowledgements}
The authors would like to thank Dan Molzahn for discussions and feedback regarding the slides. Line Roald gratefully acknowledges financial support from the Center of Non-Linear Studies (CNLS) at Los Alamos National Laboratory.

%\section*{References}
%\clearpage

LA-UR-18-26601


\begin{thebibliography}{10}

\bibitem{1512.07315}
Daniel Bienstock and Abhinav Verma.
\newblock Strong np-hardness of ac power flows feasibility.
\newblock \url{https://arxiv.org/abs/1512.07315}, 2015.

\bibitem{6581918}
W.A Bukhsh, A~Grothey, K.IM. McKinnon, and P.A Trodden.
\newblock Local solutions of the optimal power flow problem.
\newblock {\em IEEE Transactions on Power Systems}, 28(4):4780--4788, Nov 2013.

\bibitem{ac_opf_origin}
M.~J. Carpentier.
\newblock Contribution a letude du dispatching economique.
\newblock Bulletin Society Francaise Electriciens, Aug. 1962.

\bibitem{nesta}
C.~Coffrin, D.~Gordon, and P.~Scott.
\newblock {NESTA, The {\sc Nicta} Energy System Test Case Archive}.
\newblock {\em CoRR}, abs/1411.0359, 2014.

\bibitem{cp_qc_fp}
C.~Coffrin, H.~Hijazi, and P.~Van~Hentenryck.
\newblock Strengthening convex relaxations with bound tightening for power
  network optimization.
\newblock In Gilles Pesant, editor, {\em Principles and Practice of Constraint
  Programming}, volume 9255 of {\em Lecture Notes in Computer Science}, pages
  39--57. Springer International Publishing, 2015.

\bibitem{qc_opf_tps}
C.~Coffrin, H.~L. Hijazi, and P.~Van Hentenryck.
\newblock The {QC} relaxation: A theoretical and computational study on optimal
  power flow.
\newblock {\em IEEE Transactions on Power Systems}, 31(4):3008--3018, July
  2016.

\bibitem{LPAC_ijoc}
C.~Coffrin and P.~Van~Hentenryck.
\newblock A linear-programming approximation of ac power flows.
\newblock {\em INFORMS Journal on Computing}, 26(4):718--734, 2014.

\bibitem{glover2011power}
J.D. Glover, M.~Sarma, and T.~Overbye.
\newblock {\em Power System Analysis \& Design, SI Version}.
\newblock Cengage Learning, 2011.

\bibitem{grainger1994power}
J.J. Grainger and W.D. Stevenson.
\newblock {\em Power system analysis}.
\newblock McGraw-Hill series in electrical and computer engineering: Power and
  energy. McGraw-Hill, 1994.

\bibitem{932273}
I.~A. Hiskens and R.~J. Davy.
\newblock Exploring the power flow solution space boundary.
\newblock {\em IEEE Transactions on Power Systems}, 16(3):389--395, Aug 2001.

\bibitem{1664986}
R.~A. Jabr.
\newblock Radial distribution load flow using conic programming.
\newblock {\em IEEE Transactions on Power Systems}, 21(3):1458--1459, Aug 2006.

\bibitem{6407493}
R.A. Jabr.
\newblock Optimization of ac transmission system planning.
\newblock {\em IEEE Transactions on Power Systems}, 28(3):2779--2787, Aug 2013.

\bibitem{9780070359581}
P.~Kundur.
\newblock {\em Power System Stability and Control}.
\newblock McGraw-Hill Professional, 1994.

\bibitem{5971792}
J.~Lavaei and S.H. Low.
\newblock Zero duality gap in optimal power flow problem.
\newblock {\em IEEE Transactions on Power Systems}, 27(1):92 --107, feb. 2012.

\bibitem{7063278}
K.~Lehmann, A.~Grastien, and P.~Van Hentenryck.
\newblock Ac-feasibility on tree networks is np-hard.
\newblock {\em IEEE Transactions on Power Systems}, 31(1):798--801, Jan 2016.

\bibitem{6120344}
B.C. Lesieutre, D.K. Molzahn, A.R. Borden, and C.L. DeMarco.
\newblock Examining the limits of the application of semidefinite programming
  to power flow problems.
\newblock In {\em 49th Annual Allerton Conference on Communication, Control,
  and Computing (Allerton), 2011}, pages 1492 --1499, sept. 2011.

\bibitem{6815671}
S.H. Low.
\newblock Convex relaxation of optimal power flow - part ii: Exactness.
\newblock {\em IEEE Transactions on Control of Network Systems}, 1(2):177--189,
  June 2014.

\bibitem{6822653}
R.~Madani, S.~Sojoudi, and J.~Lavaei.
\newblock Convex relaxation for optimal power flow problem: Mesh networks.
\newblock {\em IEEE Transactions on Power Systems}, 30(1):199--211, Jan 2015.

\bibitem{7879340}
D.~K. Molzahn.
\newblock Computing the feasible spaces of optimal power flow problems.
\newblock {\em IEEE Transactions on Power Systems}, 32(6):4752--4763, Nov 2017.

\bibitem{forthcoming}
D.~K. Molzahn and I.~Hiskens.
\newblock A survey of relaxations and approximations of the power flow
  equations.
\newblock {\em Foundations and Trends in Electric Energy Systems}, forthcoming.

\bibitem{7447798}
D.~K. Molzahn, C.~Josz, I.~A. Hiskens, and P.~Panciatici.
\newblock A laplacian-based approach for finding near globally optimal
  solutions to opf problems.
\newblock {\em IEEE Transactions on Power Systems}, 32(1):305--315, Jan 2017.

\bibitem{molzahn2018}
D.~K. Molzahn and L.~A. Roald.
\newblock {Towards an AC Optimal Power Flow Algorithm with Robust Feasibility
  Guarantees}.
\newblock In {\em {Power Systems Computation Conference (PSCC)}}, Dublin,
  Ireland, June 2018.

\bibitem{sipser2013introduction}
Michael Sipser.
\newblock {\em Introduction to the theory of computation}.
\newblock Cengage Learning, Boston, MA, 2013.

\bibitem{6345272}
S.~Sojoudi and J.~Lavaei.
\newblock Physics of power networks makes hard optimization problems easy to
  solve.
\newblock In {\em Power and Energy Society General Meeting, 2012 IEEE}, pages
  1--8, July 2012.

\bibitem{pglib_opf}
{The IEEE PES Task Force on Benchmarks for Validation of Emerging Power System
  Algorithms}.
\newblock {PGLib Optimal Power Flow Benchmarks}.
\newblock Published online at
  \url{https://github.com/power-grid-lib/pglib-opf}.
\newblock Accessed: October 4, 2017.

\end{thebibliography}
\end{document}